\documentclass[11pt]{amsart}

\usepackage{amscd}
\usepackage{amsmath}
\usepackage{graphicx}
\usepackage{amsfonts}
\usepackage{amssymb}
\begin{document}

\title[Strong $k$-commutativity preserving maps]{Strong $k$-commutativity preserving maps on
2$\times$2 matrices }

\author{Meiyun Liu, Jinchuan Hou}
\address{Department of Mathematics, Taiyuan University of Technology,
Taiyuan 030024, P. R. of China} \email{liumeiyunmath@163.com;
jinchuanhou@aliyun.com}

\thanks{{\it 2010 Mathematics Subject Classification.}
47B49; 47B47.}
\thanks{{\it Key words and phrases.}
$k$-commutators,
 the algebra of 2$\times$2 real or complex matrices, preservers}

\thanks{This work is partially supported by  Natural
Science Foundation of China (11171249, 11271217).}

\begin{abstract}

Let ${\mathcal M}_2(\mathbb F)$ be the algebra of 2$\times$2
matrices over the real or  complex field $\mathbb F$. For a given
positive integer $k\geq 1$, the $k$-commutator of $A$ and $B$ is
defined by $[A,B]_k=[[A,B]_{k-1},B]$ with $[A,B]_0=A$ and
$[A,B]_1=[A,B]=AB-BA$.  The main result is shown that a map $\Phi:
{\mathcal M}_2(\mathbb F)\to {\mathcal M}_2(\mathbb F)$ with range
containing all rank one matrices satisfies that $[\Phi(A),\Phi(B)]_k
= [A,B]_k $ for all $A, B\in{\mathcal M}_2(\mathbb F)$ if and only
if there exist
 a functional $h :{\mathcal M}_2(\mathbb F) \rightarrow {\mathbb F}$ and a scalar
$\lambda \in{\mathbb F}$ with $\lambda^{k+1} = 1$ such that $\Phi(A)
= \lambda A + h(A)I$ for all $A \in{\mathcal M}_2(\mathbb F)$.

\end{abstract}
\maketitle

\section{Introduction}
Let ${\mathcal R}$ be a ring (or an algebra over a field ${\mathbb
F}$). For a positive integer $k\geq 1$,  recall that the
$k$-commutator of elements $A, B\in{\mathcal R}$ is defined by
$[A,B]_k = [[A,B]_{k-1},B]$ with $[A,B]_0 = A$ and $[A,B]_1 =[A,B]=
AB - BA$ being the commutator (Ref. \cite{LC}). A map
$\Phi:{\mathcal R} \rightarrow {\mathcal R}$ is said to be
commutativity preserving if $[\Phi(A), \Phi(B)] = 0$ whenever $[A,
B] = 0$ for all $A$, $B\in{\mathcal R}$; be strong $k$-commutativity
preserving if $ [\Phi(A),\Phi(B)]_k=[A,B]_k $ for all
$A,B\in{\mathcal R}$.
  Strong $1$-commutativity preserving maps are usually called strong
commutativity preserving maps. Clearly, strong commutativity
preserving maps must be commutativity preserving maps, but the
inverse is not true.

 The problem
of characterizing commutativity preserving maps have been studied
 intensively (Ref. \cite{B, BC, MS, Z} and the references therein). But in general the  commutativity preserving maps may have
uncontrollable structures. The conception of strong commutativity
preserving maps was introduced  by Bell and Dail in \cite{BD}. Then
the strong commutativity preserving additive or linear maps on
various algebraic sets were studied (Ref.  \cite{BD, BM, DA, LL}).
The results obtained there reveal that the additive and linear
strong commutativity preserving maps have nice structures.   For
general strong commutativity preserving maps without linearity or
additivity assumption,  Qi and Hou in \cite{QXH} proved that, if
$\mathcal R$ is a prime unital ring containing a nontrivial
idempotent element, then every surjective strong commutativity
preserving map $\Phi:{\mathcal R} \rightarrow {\mathcal R}$  has a
nice form, too. In fact, such map has the form $\Phi(A) = \lambda A
+ f (A)$ for all $A \in{\mathcal R}$, where $\lambda \in\{-1, 1\}$
and $f$ is a central valued map, that is, a map from ${\mathcal R}$
into its center ${\mathcal Z}_{\mathcal R}$. Liu in \cite{L}
obtained that a surjective strong commutativity preserving map
$\Phi$ on a von Neumann algebras ${\mathcal A}$ without central
summands of type $I_1$ has the form $\Phi(A) = ZA + f(A)$ for all $A
\in{\mathcal A}$, where $Z\in{\mathcal Z_{\mathcal A}}$ with $Z^2 =
I$ and $f$ is a central valued map.

It seems that the study of the problem of characterizing strong
$k$-commutativity preserving maps for $k>1$ was started by
 \cite{Q}, in which the form of strong $2$-commutativity preserving map is given  on a
 unital prime ring containing a nontrivial idempotent.
Recently, the strong $3$-commutativity preserving map was
characterized on standard operator algebra ${\mathcal A}$ in
${\mathcal B}(X)$ in \cite{LH}, where $X$ is a Banach space of
dimension $\geq 2$ over the real or complex field ${\mathbb F}$. The
result in \cite{LH} revealed  that, if $\Phi$ is a surjective map on
${\mathcal A}$, then $\Phi$ is strong $3$-commutativity preserving
if and only if there exist
 a functional $h :{\mathcal A} \rightarrow {\mathbb F}$ and a scalar
$\lambda \in{\mathbb F}$ with $\lambda^4 = 1$ such that $\Phi(A) =
\lambda A + h(A)I$ for all $A \in{\mathcal A}$.
 It is natural to raise the problem how to characterize the strong $k$-commutativity
preserving maps  on rings for any positive integer $k$. Notice that,
$[A,B]_k=\sum_{i=0}^k(-1)^i{\rm C}_k^iB^iAB^{k-i}$, where, ${\rm
C}_k^i=\frac{k(k-1)\cdots(k-i+1)}{i!}$. Thus, with $k$ increasing,
the problem of characterizing strong $k$-commutativity preserving
maps becomes much more difficult. The purpose of this paper is  to
answer the above problem for the case when   the maps act on the
algebra of 2$\times$2 matrices over the real or  complex field
$\mathbb F$.

The following is our main results.

\textbf {Theorem 1.} {\it Let ${\mathcal M}_2(\mathbb F)$ be a
 algebra of 2$\times$2
matrices over the real or  complex field $\mathbb F$ and $k\geq 1$
be an integer. Assume that $\Phi:{\mathcal M}_2(\mathbb
F)\rightarrow {\mathcal M}_2(\mathbb F)$ is a map with range
containing all rank one matrices. Then $\Phi$ is strong
$k$-commutativity preserving if and only if there exist a functional
$h :{\mathcal M}_2(\mathbb F)\rightarrow {\mathbb F}$ and a scalar
$\lambda \in{\mathbb F}$ with $\lambda^{k+1} = 1$ such that $\Phi(A)
= \lambda A + h(A)I$ for all $A \in{\mathcal M}_2(\mathbb F)$.}

\section{Proof of the main result}

 Before proving Theorem 1, we  give several  lemmas.

The first lemma is obvious by the main result in \cite{FS}.

 \textbf{Lemma 2.1.}  {\it Let ${\mathcal M}_2(\mathbb F)$ be a
 algebra of 2$\times$2
matrices over the real or  complex field $\mathbb F$, and let
$\{A_i, B_i\}_{i=1}^n$,
  $\{C_j,D_j\}_{j=1}^m\subset{\mathcal M}_2(\mathbb F)$,
such that $\sum{_{i=1}^n}
 A_{i}T
B_{i} = \sum{_{i=1}^n} C_{j}TD_{j}$ for all rank one $ T\in{\mathcal
M}_2(\mathbb F)$. If $A_{1},\ldots,A_{n}$ are linearly independent,
then each $B_{i}$ is a linear combination of $D_{1},\ldots,D_{m}$.
Similarly, if $B_{1},\ldots,B_{m}$ are linearly independent, then
each $A_{i}$ is a linear combination of $C_{1},\ldots,C_{m}$.}

\textbf{Lemma 2.2.} {\it Let positive integer $k\geq 1$ and
$Z\in{\mathcal M}_2(\mathbb F)$. Assume $[Z,A]_k = 0$ for any rank
one idempotent matrix
 $A\in{\mathcal M}_2(\mathbb F)$ , then there exists a scalar
$\lambda\in{\mathbb F}$ such that $Z=\lambda I$, where $I$ is the
identity matrix.}

 \textbf{Proof.}  Every rank one matrix $A$ can be written as
 $A=xf^*$ for some $x,f\in{\mathbb F}^{(2)}$, where $B^*$ stands for the conjugate transpose of the matrix $B$.
Clearly,  $A=xf^*$ is an idempotent if and only if $\langle
x,f\rangle=f^*x=1$. Thus by the assumption, $$[Z,xf^*]_k =
\sum_{i=1}^k(-1)^i{\rm C}_{k}^i (xf^*)^iZ(xf^*)^{k-i}=0$$ for any
$x,f\in {\mathbb F}^{(2)}$ with $\langle x,f\rangle=1$. As
$(xf^*)^0=I$ and $(xf^*)^i=xf^*$ when $i\geq 1$, one sees that, if
$k$ is odd, then
$$[Z,xf^*]_k =Zxf^*-x(Z^*f)^*=0$$ for all
$x,f\in {\mathbb F}^{(2)}$ with $\langle x,f\rangle=1$; if $k$ is
even, then
$$[Z,xf^*]_k =Zxf^*-2
 (xf^*)Z(xf^*)+xf^*Z= Zxf^*-2\langle Zx,f\rangle xf^*+x(Z^*f)^*=0$$ for all $x,f\in {\mathbb
 F}^{(2)}$ with $\langle x,f\rangle=1$. Above two identities imply
 that, in any case, $Z$ has the property that $Zx=\lambda_x x$ holds
 for all $x\in{\mathbb F}^{(2)}$. It follows that, there must be a
 scalar $\lambda$ such that $Z=\lambda I$. \hfill$\Box$

Denote by $\mathcal N({\mathcal M}_2(\mathbb F))$ the set of all
nilpotent element in ${\mathcal M}_2(\mathbb F)$.

 \textbf{Lemma 2.3.} {\it Let $k\geq 3$ be a  positive integer and $S \in{\mathcal M}_2(\mathbb F)$. Then $[A,S]_k = 0$
holds for any rank one matrix $A \in{\mathcal M}_2(\mathbb F)$ if
and only if there exists a scalar $\lambda \in{\mathbb F}$ and an
element $N \in{\mathcal N({\mathcal M}_2(\mathbb F)})$ such that $S
= \lambda I + N$.}

 \textbf{ Proof.} To check the ``if" part, assume that $S = \lambda I +
N$ with $N^2 = 0$. It is easily seen that
$$[A, S]_k = [A, \lambda I + N]_k = [A, N]_k =\sum_{i=0}^k(-1)^i{\rm C}_k^iN^iAN^{k-i}= 0$$
for any $ A\in{\mathcal M}_2(\mathbb F).$

Next we check the `` only if " part.

\textbf{Case 1.} ${\mathbb F} = {\mathbb C}$.

Since $S\in{\mathcal M}_2(\mathbb C)$,  there exists a quadratic
polynomial $P(t) = (t-\alpha_1)(t-\alpha_2)$, such that $P(S) = (S
-\alpha_1 I)(S -\alpha_2 I) = 0$.

 If $\alpha_1 \neq \alpha_2$, taking eigenvectors $x_i$, $f_i\in{\mathbb C}^{(2)}$ of $S$ and $S^*$ with
respect to  $\alpha_i$ such that $Sx_i=\alpha_ix_i$ and
$S^*f_i=\bar{\alpha}_if_i$, $i=1,2$.

Taking $A=x_1f_2^*$, then
$$\begin{array}{rl}
0=&[A, S]_k = [x_1f_2^*, S]_k
=\sum_{i=0}^k(-1)^i{\rm C}_k^iS^i(x_1 f_2^*)S^{k-i}\\
=&\sum_{i=0}^k(-1)^i{\rm C}_k^iS^ix_1({S^*}^{k-i}f_2)^*
=\sum_{i=0}^k(-1)^i{\rm
C}_k^i\alpha_1^i\bar{\alpha}_2^{k-i}x_1f_2^*\\
=&\sum{_{i=0}^k}(-1)^i{\rm C}_k^i\alpha_1^i \alpha_2^{k-i}x_1 f_2^*=
(\alpha_2 - \alpha_1)^kx_1 f_2^*\not=0,
\end{array}$$
a contradiction.  So we must have $\alpha_1 = \alpha_2$, and then
$P(S) = (S -\alpha_1 I)^2 = 0$. Let $N=S -\alpha_1 I$, $\lambda =
\alpha_1$, then $S=\lambda I + N$ with $N^2=0$.

\textbf{Case 2.}  ${\mathbb F} = {\mathbb R}$.

 Since $S\in{\mathcal
M}_2(\mathbb R)$,   there exists quadratic polynomial $P(t) =
(t-\alpha_1)(t-\alpha_2)$ or $P(t) = (t + \alpha)^2 + \beta^2$,
where $\alpha,\beta\in{\mathbb R}$ with $\beta\not=0$, such that
$P(S) = 0$.

If $P(t) = (t-\alpha_1)(t-\alpha_2)$, by a similar argument as the
proof for the case when ${\mathbb F} = {\mathbb C}$, we know that
$S$ is of the form $S = \lambda I + N$ with $N$ nilpotent.

We claim that the case $P(t) = (t + \alpha)^2 + \beta^2$ does not
happen. Assume, on the contrary,  $P(t) = (t + \alpha)^2 + \beta^2$;
then $(S + \alpha I)^2 = -\beta^2I$, which implies that $S$ is not
of the form $\lambda I$ for any $\lambda\in{\mathbb R}$. Let
$S^\prime = S + \alpha I$; then ${S^\prime}^2 = -\beta^2I$. For any
rank one $A\in {\mathcal M}_2({\mathbb R})$, we have
$$\sum_{i=0}^k(-1)^i{\rm C}_k^i{S^\prime}^iA{S^\prime}^{k-i}=[A, S^\prime]_k = [A, S + \alpha I]_k
=[A, S]_k =0.$$

If $k$ is odd, one gets
$$\begin{array}{rl}
 [A, S^\prime]_k=&{\rm C}_k^0A(-\beta^2)^{\frac{k-1}{2}}{S^\prime}
-{\rm C}_k^1{S^\prime}A(-\beta^2)^{\frac{k-1}{2}}I + \ldots + {\rm
C}_k^{k-1}(-\beta^2)^{\frac{k-1}{2}}A{S^\prime}
-{\rm C}_k^k(-\beta^2)^{\frac{k-1}{2}}{S^\prime}A\\
=&(-\beta^2)^{\frac{k-1}{2}}({\rm C}_k^0A{S^\prime}-{\rm
C}_k^1{S^\prime}A + {\rm C}_k^2A{S^\prime}- {\rm C}_k^3{S^\prime}A +
\ldots +
{\rm C}_k^{k-1}A{S^\prime}- {\rm C}_k^k{S^\prime}A)\\
=&(-\beta^2)^{\frac{k-1}{2}}({\rm C}_k^0 + {\rm C}_k^1 + \ldots +
{\rm C}_k^{\frac{k-1}{2}})(A{S^\prime}- {S^\prime}A)\\
=&(-\beta^2)^{\frac{k-1}{2}}2^{k-1}(A{S^\prime}- {S^\prime}A)\\
=&(-\beta^2)^{\frac{k-1}{2}}2^{k-1}(AS- SA)=0,
\end{array}$$
which implies that $AS=SA$ holds for all rank one matrix $A\in
{\mathcal M}_2({\mathbb R})$. So there exists scalar
$\lambda\in{\mathbb R}$, such that $S = \lambda I$, a contradiction.

 If $k$ is even, we have
$$\begin{array}{rl}
[A, S^\prime]_k =& {\rm C}_k^0A(-\beta^2)^{\frac{k}{2}} -{\rm
C}_k^1{S^\prime}A(-\beta^2)^{\frac{k-2}{2}}S^\prime + \ldots - {\rm
C}_k^{k-1}(-\beta^2)^{\frac{k-2}{2}}{S^\prime}A{S^\prime}
+ {\rm C}_k^k(-\beta^2)^{\frac{k}{2}}A\\
=&(-\beta^2)^{\frac{k-2}{2}}[-\beta^2({\rm C}_k^0 + {\rm C}_k^2 +
\ldots +
{\rm C}_k^k)A - ({\rm C}_k^1 +{\rm C}_k^3 + \ldots + {\rm C}_k^{k-1}){S^\prime}A{S^\prime}]\\
=&(-\beta^2)^{\frac{k-2}{2}}2^{k-1}(-\beta^2A-{S^\prime}A{S^\prime})
= 0.
\end{array}$$
So, we have $\beta^2A + {S^\prime}A{S^\prime} = 0$ for any rank one
matrix. Then, applying Lemma 2.1, there exists scalar
$\lambda_1\in{\mathbb R}$ such that $S+ \alpha I=S^\prime =
\lambda_1I$. So $S = \lambda I$ with $\lambda=\lambda_1-\alpha$,
again getting a contradiction.

 The proof of the Lemma 2.3 is completed.
\hfill$\Box$

Now we are at the position to give our proof of the main theorem.

\textbf {Proof of Theorem 1.} The `` if " part of the theorem is
obvious. \if As ${\mathcal M}_2({\mathbb F})$ is a prime algebra,
 the proof of the theorem 1 had been finished in \cite{QXH} and
\cite{Q} for the cases when $k=1,2$. Thus, to check the ``only if"
part, it suffices to check the case when $k\geq 3$.\fi

In the sequel, we always assume that $k\geq 1$ and $\Phi:{\mathcal
M}_2(\mathbb F) \rightarrow {\mathcal M}_2(\mathbb F)$ is a   strong
$k$-commutativity preserving map with range containing all rank one
matrices. We will check the `` only if " part by several steps.

\textbf {Step 1.} For any $A, B\in{\mathcal M}_2(\mathbb F)$, there
exists a scalar $\lambda_{A,B} \in{\mathbb F}$ such that $\Phi(A +
B) = \Phi(A) + \Phi(B) + \lambda_{A,B} I$.

 For any $A, B$ and $T\in{\mathcal M}_2(\mathbb F)$, we have
 \begin{align*}
 &[\Phi(A + B) - \Phi(A) - \Phi(B), \Phi(T)]_k\\
 &= [\Phi(A + B), \Phi(T)]_k - [\Phi(A), \Phi(T)]_k - [\Phi(B),
 \Phi(T)]_k\\
&= [A  + B, T]_k - [A, T]_k - [B, T]_k = 0.
\end{align*}
Since the range of $\Phi$ contains all rank one matrices, we can
apply Lemma 2.2 to ensure that the above equation implies $\Phi(A +
B) - \Phi(A) - \Phi(B)\in\{\lambda I:\lambda\in{\mathbb F}\}$. So
the assertion in Step 1 is true.

\textbf {Step 2.} $\Phi({\mathbb F}I) = {\mathbb F}I \cap {\rm ran}
(\Phi)$ and $\Phi({\mathbb F}I + {\mathcal N({\mathcal M}_2(\mathbb
F))}) = ({\mathbb F}I + {\mathcal N}({\mathcal M}_2(\mathbb F)))\cap
{\rm ran} (\Phi)$.

 $\Phi({\mathbb F}I) = {\mathbb F}I\cap {\rm ran}
(\Phi) $ is obvious by Lemma 2.2 and the assumption on
 the range of $\Phi$.

If $N\in {\mathcal N}({\mathcal M}_2(\mathbb F))$ and
$\lambda\in{\mathbb F}$, then, for any $A\in{\mathcal M}_2(\mathbb
F)$, we have
$$[\Phi(A), \Phi(\lambda I + N)]_k = [A, \lambda I + N]_k
= [A, N]_k = 0.$$
 By the assumption on the range of $\Phi$ and Lemma 2.3, we
see that $\Phi(\lambda I + N) \in {\mathbb F}I + {\mathcal
N}({\mathcal M}_2(\mathbb F))$.

On the other hand, if $\Phi(B) = \lambda I + N$ for some scalar
$\lambda$ and $N\in{\mathcal N}({\mathcal M}_2(\mathbb F))$, then
$$[A,B]_k = [\Phi(A),\Phi(B)]_k = [\Phi(A), \lambda I + N]_k = [\Phi(A),
N]_k = 0.$$ for all $A\in{\mathcal M}_2(\mathbb F)$. By Lemma 2.3
again, we see that $B \in {\mathbb F} I + {\mathcal N}({\mathcal
M}_2(\mathbb F))$.

Hence $\Phi({\mathbb F} I + {\mathcal N}({\mathcal M}_2(\mathbb F)))
=( {\mathbb F}I + {\mathcal N}({\mathcal M}_2(\mathbb F)))\cap {\rm
ran} (\Phi)$,
 completing the proof of the Step 2.

Let $$E_{11} = \left(
\begin{array}{cc}
  1 &0 \\
  0&0
\end{array}
\right), \ \ E_{22} = \left(
\begin{array}{cc}
  0&0 \\
  0&1
\end{array}
\right),\ \  E_{12} = \left(
\begin{array}{cc}
  0&1 \\
  0&0
\end{array}
\right)\ \mbox{\rm and}\  E_{21} = \left(
\begin{array}{cc}
  0&0 \\
 1&0
\end{array}
\right).$$
 Then every $A\in {\mathcal M}_2(\mathbb F)$ can be written as $A =
a_{11}E_{11} + a_{12}E_{12} + a_{21}E_{21} + a_{22}E_{22}$, where
${a}_{ij} \in{\mathbb F}$.

\textbf{Step 3.} There exist two scalars $\lambda$,
$\mu_1\in{\mathbb F}$ with $\lambda\not=0$, such that
 $\Phi(E_{11}) = \lambda E_{11} + \mu_1 I$.

We prove the assertion in Step 3 by three claims.

\textbf{Claim 3.1.} For any $A\in{\mathcal M}_2(\mathbb F)$, we have
$[A, \Phi(E_{11})]_k\in {\mathbb F}E_{12} + {\mathbb F}E_{21}$.

Note that, for any $A\in{\mathcal M}_2(\mathbb F)$, we have
$$[A, E_{11}]_2 = a_{12}E_{12} + a_{21}E_{21}\in{\mathbb
F}E_{12} + {\mathbb F}E_{21}.$$ Also notice  that $[A, Q] = [A,
Q]_3$ holds for any $A, Q\in {\mathcal A}$ with $Q$ an idempotent.
So we have  $[A, E_{11}]_k = [A, E_{11}]_{k+2} = [[A, E_{11}]_k,
E_{11}]_2$, which implies that $$[\Phi(A), \Phi(E_{11})]_k =
[[\Phi(A), \Phi(E_{11})]_k, E_{11}]_2$$ holds for any $A\in{\mathcal
M}_2(\mathbb F)$.  Since the range of $\Phi$ contains all rank one
matrices and every matrix in ${\mathcal M}_2({\mathbb F})$ is a sum
of rank one matrices, one gets
$$[A, \Phi(E_{11})]_k = [[A, \Phi(E_{11})]_k, E_{11}]_2\ \ \mbox{\rm for all}\ A\in{\mathcal M}_2(\mathbb F).$$
Therefore, for any $A$ we have $$[A, \Phi(E_{11})]_k\in{\mathbb
F}E_{12} + {\mathbb F}E_{21}.$$

 \textbf{Claim 3.2.}  The result of Step 3 is true when ${\mathbb F} = {\mathbb
C}$.

Assume ${\mathbb F} = {\mathbb C}$. Since $\Phi(E_{11})\in{\mathcal
M}_2(\mathbb C)$, there exists a polynomial
$p(t)=(t-\alpha)(t-\beta)$ such that $p(\Phi(E_{11}))=0$.

If $\alpha = \beta$, then $\Phi(E_{11})$ can be written as
$\Phi(E_{11})=\alpha I+N$ with $N^2=0$. But, by Step 2, this entails
that $E_{11}\in{\mathbb F}I+{\mathcal N}({\mathbb F})$, a
contradiction. So $\alpha\neq\beta$. Let $x$ be an eigenvector of
$\Phi(E_{11})$ with respect to  $\alpha$ and $f$ be an eigenvector
of $\Phi(E_{11})^*$ with respect to $\bar{\beta}$;  then
$\Phi(E_{11})x =\alpha x$ and $\Phi(E_{11})^*f = \bar{\beta} f$.
Taking $B =x  f^*$ and applying Claim 3.1 give
$$\begin{array}{rl}
[B, \Phi(E_{11})]_k &= [xf^*, \Phi(E_{11})]_k\\
&=\sum_{i=0}^k(-1)^i{\rm C}_k^i\Phi(E_{11})^i(x f^*)\Phi(E_{11})^{k-i}\\
&=\sum_{i=0}^k(-1)^i{\rm C}_k^i\Phi(E_{11})^ix({\Phi(E_{11})^*}^{k-i}f)^*\\
&=\sum_{i=0}^k(-1)^i{\rm C}_k^i\alpha^ix (\overline{\beta}^{k-i}f)^*\\
&=\sum_{i=0}^k(-1)^i{\rm C}_k^i\alpha^i \beta^{k-i}x  f^*\\
&=(\beta-\alpha)^kx  f^*\in {\mathbb C}E_{12}+{\mathbb C}E_{21}.
\end{array}$$
Multiplying $E_{11}$ from both sides of the above equation, one gets
$E_{11}(x f^*)E_{11} = E_{11}x ( E_{11}f)^* = 0$; similarly,
multiplying $E_{22}$ from both sides of the above equation, one gets
$E_{22}(x f^*)E_{22} = E_{22}x( E_{22}f)^* = 0$. Hence, we have
$$\begin{cases}
E_{11}x=0\ \ \mbox{\rm or}\ \ E_{11}f = 0,\\
E_{22}x=0\ \ \mbox{\rm or}\ \ E_{22}f = 0.
\end{cases}
\eqno(2.1)$$ Let  $y$, $g$  be respectively  eigenvectors of
$\Phi(E_{11})$ and $\Phi(E_{11})^*$ with respect to $\beta$ and
$\alpha$, a similar argument as above gives
$$\begin{cases}
E_{11}y=0\ \ \mbox{\rm or}\ \ E_{11}g = 0,\\
E_{22}y=0\ \ \mbox{\rm or}\ \ E_{22}g = 0.
\end{cases}
\eqno(2.2)$$
 By  Eq.$(2.1)\sim(2.2)$ and linear independence  of $x$ and $y$,
 we get
$$x = \left(\begin{array}{c}
  1 \\
  0
\end{array}\right),\ \  y = \left(\begin{array}{c}
  0 \\
  1
\end{array}\right)\ \ \mbox{\rm or}\
\ \ x = \left(\begin{array}{c}
  0\\
  1
\end{array}\right),\ \ y = \left(\begin{array}{c}
  1 \\
  0
\end{array}\right).$$
So $\Phi(E_{11})$ is diagonal, say $\Phi(E_{11}) = \left(
\begin{array}{cc}
  \alpha &0 \\
  0&\beta
\end{array}
\right)$. Let $\lambda
 = \alpha-\beta$ and $\mu_1 = \beta$; then
  $$\Phi(E_{11}) = \lambda E_{11} + \mu_1 I.$$
Finally, by Step 2, we see that $\lambda \not=0$, as desired.

\textbf{Claim 3.3.}  The assertion of Step 3 is also true for
 the case ${\mathbb F} = {\mathbb R}$.

Since $\Phi(E_{11})\in{\mathcal M}_2(\mathbb R)$,  there exists a
quadratic polynomial $p(t) = (t-\alpha_1)(t-\alpha_2)$ or $P(t) = (t
+ \alpha)^2 + \beta^2$,
 where $\beta\neq0$, such that $p(\Phi(E_{11})) = 0$.

If $P(t) = (t-\alpha_1)(t-\alpha_2)$, then $p(\Phi(E_{11})) =
(\Phi(E_{11})-\alpha_1I)(\Phi(E_{11})-\alpha_2I) = 0$. Then a
similar argument to the case when ${\mathbb F} = {\mathbb C}$ shows
that there exist $\lambda, \mu_1\in{\mathbb R}$ such as
$\Phi(E_{11}) = \lambda E_{11} + \mu_1 I$.

Next we show that the case that $p(t)$ has the form of $p(t) = (t +
\alpha)^2 + \beta^2$ never occurs. If, on the contrary,  $p(t) = (t
+ \alpha)^2 + \beta^2$, then we have $(\Phi(E_{11}) + \alpha I)^2 =
-\beta^2I$.  Write $\Phi(E_{11}) = s_{11}E_{11} + s_{12}E_{12} +
s_{21}E_{21} + s_{22}E_{22}$ with $s_{ij} \in{\mathbb R}$ and let
$\Phi(E_{11})^\prime = \Phi(E_{11}) + \alpha I$; then
${\Phi(E_{11})^\prime}^2 = -\beta^2I$. For any $B\in{\mathcal
M}_2(\mathbb R)$, we have
$$[B, \Phi(E_{11})]_k = [B, \Phi(E_{11}) + \alpha I]_k = [B, \Phi(E_{11})^\prime]_k
=\sum_{i=0}^k(-1)^i{\rm
C}_k^i{\Phi(E_{11})^\prime}^iB{\Phi(E_{11})^\prime}^{k-i}.$$

If $k$ is odd, by Claim 3.1, we have
$$\begin{array}{rl}
&[B, \Phi(E_{11})]_k\\
 = &{\rm C}_k^0B(-\beta^2)^{\frac{k-1}{2}}\Phi(E_{11})^\prime
-{\rm C}_k^1\Phi(E_{11})^\prime B(-\beta^2)^{\frac{k-1}{2}}I +
\ldots
-{\rm C}_k^k(-\beta^2)^{\frac{k-1}{2}}\Phi(E_{11})^\prime B\\
=&(-\beta^2)^{\frac{k-1}{2}}({\rm C}_k^0B\Phi(E_{11})^\prime-{\rm
C}_k^1\Phi(E_{11})^\prime B + \ldots +
{\rm C}_k^{k-1}B\Phi(E_{11})^\prime - {\rm C}_k^k\Phi(E_{11})^\prime B)\\
=&(-\beta^2)^{\frac{k-1}{2}}({\rm C}_k^0 + {\rm C}_k^1 + \ldots +
{\rm C}_k^{\frac{k-1}{2}})(B \Phi(E_{11})^\prime - \Phi(E_{11})^\prime B)\\
=&(-\beta^2)^{\frac{k-1}{2}}2^{k-1}(B \Phi(E_{11})^\prime -
 \Phi(E_{11})^\prime B)\\
=&(-\beta^2)^{\frac{k-1}{2}}2^{k-1}(B\Phi(E_{11})-
\Phi(E_{11})B)\in{\mathbb R}E_{12} + {\mathbb R}E_{21}.
\end{array} \eqno(2.3)$$
 Taking $B=E_{21}$ in Eq.(2.3) gives
 $$ [E_{21}, \Phi(E_{11})]_k =
(-\beta^2)^{\frac{k-1}{2}}2^{k-1}[
 (s_{11}- s_{22})E_{21} + s_{12}E_{22} -  s_{12}E_{11}]\in {\mathbb R}E_{12} +
{\mathbb R}E_{21}.$$  So $s_{12} = 0$.  Taking $B=E_{12}$ in
Eq.(2.3), one obtains
$$ [E_{12}, \Phi(E_{11})]_k =
 (-\beta^2)^{\frac{k-1}{2}}2^{k-1}[s_{21}E_{11} + (s_{22} -  s_{11})E_{12}- s_{21}E_{22}] \in {\mathbb R}E_{12} +
{\mathbb R}E_{21},$$ which forces $s_{21} = 0$. Hence we have
$\Phi(E_{11})=s_{11}E_{11} + s_{22}E_{22} $, which contradicts to
the fact that $(\Phi(E_{11})+\alpha I)^2+\beta^2 I=0$.

 If $k$ is even, by applying Claim 3.1, we have
$$\begin{array}{rl}
&[B, \Phi(E_{11})]_k=[B, \Phi(E_{11})^\prime]_k \\=&{\rm
C}_k^0B(-\beta^2)^{\frac{k}{2}}-{\rm C}_k^1\Phi(E_{11})^\prime
B(-\beta^2)^{\frac{k-2}{2}}I + \ldots
+{\rm C}_k^k(-\beta^2)^{\frac{k}{2}}B\\
=&(-\beta^2)^{\frac{k-2}{2}}[-\beta^2({\rm C}_k^0+{\rm C}_k^2 +
\ldots +
{\rm C}_k^k)B -({\rm C}_k^1+{\rm C}_k^3+\ldots +{\rm C}_k^{k-1})\Phi(E_{11})^\prime B \Phi(E_{11})^\prime]\\
=&(-\beta^2)^{\frac{k-2}{2}}2^{k-1}(-\beta^2B-\Phi(E_{11})^\prime
B\Phi(E_{11})^\prime)\in{\mathbb R}E_{12} + {\mathbb R}E_{21}.
\end{array}\eqno(2.4)$$
Write $\Phi(E_{11})^\prime = s_{11}^\prime E_{11} + s_{12} ^\prime
E_{12} + s_{21} ^\prime E_{21} + s_{22}^\prime E_{22}$ with
$s_{ij}^\prime \in{\mathbb R}$.
 Take $B=E_{11}$ in Eq.(2.4); then
$$\begin{array}{rl}
&[E_{11}, \Phi(E_{11})^\prime]_k\\ =&
(-\beta^2)^{\frac{k-2}{2}}2^{k-1}[(-\beta^2
-{s_{11}^\prime}^2)E_{11}- s_{11}^\prime s_{12}^\prime E_{12} -
s_{21}^\prime s_{11}^\prime E_{21} - s_{21}^\prime s_{12}^\prime
E_{22}]\\ \in & {\mathbb R}E_{12} + {\mathbb R}E_{21}, \end{array}$$
which implies that ${s_{11}^\prime}^2 = -\beta^2$, a contradiction.

\if false
$$s_{21}^\prime s_{12}^\prime = 0\ \ \mbox{\rm and}\ \ {s_{11}^\prime}^2 = -\beta^2.\eqno(2.5)$$
Take $B=E_{22}$ in Eq.(2.4); then
$$\begin{array}{rl} &[E_{22},
\Phi(E_{11})^\prime]_k \\=&
(-\beta^2)^{\frac{k-2}{2}}2^{k-1}[(-\beta^2-{s_{22}^\prime}^2)E_{22}-
s_{12}^\prime s_{22}^\prime E_{12} - s_{22}^\prime s_{21}^\prime
E_{21} - s_{12}^\prime s_{21}^\prime E_{11}]\\ \in & {\mathbb
R}E_{12} + {\mathbb R}E_{21}, \end{array}$$ which entails that
$$s_{12}^\prime s_{21}^\prime = 0\ \ \mbox{\rm or}\ \ {s_{22}^\prime}^2 = -\beta^2.\eqno(2.6)$$
Take $B=E_{12}$ in Eq.(2.4); then
$$\begin{array}{rl} &[E_{12},\Phi(E_{11})^\prime]_k \\=& (-\beta^2)^{\frac{k-2}{2}}2^{k-1}[-s_{11}^\prime s_{21}^\prime E_{11} + (-\beta^2 -
s_{11}^\prime s_{22}^\prime )E_{12} - {s_{21}^\prime}^2 E_{21} -
s_{21}^\prime s_{22}^\prime E_{22}]\\ \in & {\mathbb R}E_{12} +
{\mathbb R}E_{21}, \end{array}$$ which forces that
$$s_{21}^\prime s_{22}^\prime = 0\ \ \mbox{\rm or}\ \ s_{11}^\prime
s_{21}^\prime = 0.\eqno(2.7)$$ Similarly, taking $B=E_{21}$ in
Eq.(2.4) gives
$$s_{11}^\prime s_{12}^\prime=s_{22}^\prime s_{12}^\prime=0.
\eqno(2.8)
$$
 By  Eq.$(2.5)\sim(2.8)$, we see $s_{21}^\prime = s_{12}^\prime = 0.$
So $s_{12}=s_{21}=0$ and $\Phi(E_{11})= s_{11}  E_{11} + s_{22}
E_{22}, $ again a contradiction.\fi

\textbf{Step 4.}  For $1\leq i\not = j\leq 2$ and for any
$a_{ij}\in\mathbb F$, there exists a scalar $\mu_{a_{ij}}\in{\mathbb
F}$, such that $\Phi(a_{ij}E_{ij}) = \lambda^{-k}a_{ij}E_{ij} +
\mu_{a_{ij}} I$.

Here, we only give the proof for the case when $(i,j) = (1,2)$. The
proof for $(i,j) = (2,1)$ is similar.

For any $a_{12}\in{\mathbb F}$, write $\Phi(a_{12}E_{12}) =
s_{11}E_{11} + s_{12}E_{12} + s_{21}E_{21} + s_{22}E_{22}$ with
$s_{ij}\in{\mathbb F}$.  Then, by Step 3, we have
$$\begin{array}{rl}
 (-1)^k a_{12}E_{12} &= [a_{12}E_{12}, E_{11}]_k = [\Phi(a_{12}E_{12}), \Phi(E_{11})]_k\\
 &= [s_{11}E_{11} + s_{12}E_{12} + s_{21}E_{21}+ s_{22}E_{22}, \lambda E_{11}]_k \\
 &= \lambda^k \sum_{i=1}^k(-1)^i{\rm C}_k^iE_{11}^i(s_{11}E_{11} + s_{12}E_{12} + s_{21}E_{21}+ s_{22}E_{22})E_{11}^{k-i}\\
 &= \lambda^k(s_{21}E_{21} +(-1)^k s_{12}E_{12}),
\end{array}$$
which entails that  $s_{21} = 0$ and $\lambda^k s_{12} = a_{12}$. So
$s_{12} = \lambda^{-k} a_{12}$ and
$$\Phi(a_{12}E_{12}) = \lambda^{-k} a_{12}E_{12} + s_{11}E_{11} + s_{22}E_{22}.$$
On the other hand, by Step 2,  $\Phi(a_{12}E_{12}) = \mu I + N$ for
some $\mu\in{\mathbb F}$ and $N\in{\mathcal N({\mathcal M}_2(\mathbb
F))}$. Write $N = \left(
\begin{array}{cc}
  r_{11} &r_{12}\\
  r_{21}&r_{22}
\end{array}
\right)$.
It follows that
$$\begin{array}{rl}
&\lambda^{-k} a_{12}E_{12} + s_{11}E_{11} + s_{22}E_{22} = \mu I + N \\
=& (\mu + r_{11})E_{11} + r_{12}E_{12} + r_{21}E_{21} + (\mu +
r_{22})E_{22}.
\end{array}$$
This gives that $r_{21} = 0$, $r_{12} = \lambda^{-k} a_{12}$. Since
$N^2 = 0$, we must have $r_{11} = r_{22} = 0$, and then $s_{11} =
s_{22} = \mu_{a_{12}}$.
 Thus one gets
 $$\Phi(a_{12}E_{12}) = \lambda^{-k} a_{12}E_{12} + \mu_{a_{12}} E_{11} + \mu_{a_{12}} E_{22} = \lambda^{-k} a_{12}E_{12} + \mu_{a_{12}}I,$$
 as desired.

\textbf{Step 5.}  $\lambda^{k+1} = 1$, and, for any
$a_{ii}\in{\mathbb F}_{ii}$, $(i\in\{1,2\})$, there exists a scalar
$\mu_{a_{ii}}\in{\mathbb F}$, such that $\Phi(a_{ii}E_{ii}) =
\lambda a_{ii}E_{ii} + \mu_{a_{ii}} I.$

Still, we only prove that the claim holds for the case $i=1$.

Take any nonzero $a_{11}\in{\mathbb F}$ and write
$\Phi(a_{11}E_{11}) = s_{11}E_{11} + s_{12}E_{12} + s_{21}E_{21}+
s_{22}E_{22}$ with $s_{ij}\in{\mathbb F}$. By Step 3, we obtain that
$$\begin{array}{rl}
 0 &= [a_{11}E_{11}, E_{11}]_k = [\Phi(a_{11}E_{11}), \Phi(E_{11})]_k\\
 &= [s_{11}E_{11} + s_{12}E_{12} + s_{21}E_{21}+ s_{22}E_{22}, \lambda E_{11} + \mu_1 I ]_k \\
 & = \lambda^k(s_{21}E_{21} +(-1)^k  s_{12}E_{12}),
\end{array}$$
which implies $s_{21} = s_{12} = 0$, and then $\Phi(a_{11}E_{11}) =
s_{11}E_{11} + s_{22}E_{22}$. By Step 4, there exists a scalar
$\mu_{E_{12}}\in{\mathbb F}$, such that $\Phi(E_{12})=\lambda^{-k}
E_{12}+\mu_{E_{12}}I$. Thus we have
$$\begin{array}{rl}
 (-1)^k a_{11}^kE_{12} =& [E_{12}, a_{11}E_{11}]_k = [\Phi(E_{12}), \Phi(a_{11}E_{11})]_k\\
 =& [\lambda^{-k} E_{12} + \mu_{E_{12}} I, s_{11}E_{11} + s_{22}E_{22} ]_k \\
 =&\lambda^{-k}\sum_{i=0}^k(-1)^i{\rm C}_k^i(s_{11}E_{11} + s_{22}E_{22})^i E_{12}(s_{11}E_{11} + s_{22}E_{22})^{k-i}\\
 =&\lambda^{-k}\sum_{i=0}^k(-1)^i{\rm C}_k^is_{11}^is_{22}^{k-i} E_{12}\\
 = & \lambda^{-k}(s_{22} - s_{11})^kE_{12}\\
 = & (-1)^k\lambda^{-k}(s_{11} - s_{22})^kE_{12}.
\end{array}$$
So
$$(s_{11} - s_{22})^k = \lambda^k a_{11}^k,$$
which entails that $$s_{11} - s_{22} = \delta \lambda a_{11}$$ for
some $\delta$ with $\delta^k = 1$.
 Thus we have
$$\Phi(a_{11}E_{11}) = s_{11}E_{11} + s_{22}E_{22} = (\delta \lambda a_{11} + s_{22})E_{11} + s_{22}E_{22}= \delta \lambda a_{11}E_{11} + s_{22} I.$$
Applying Step 1 and Step 4, it is easily checked that
$$\begin{array}{rl}
 &a_{11}[E_{11}, E_{12} + E_{21}]_k\\
 =& [a_{11}E_{11}, E_{12} + E_{21}]_k = [\Phi(a_{11}E_{11}), \Phi(E_{12} +
 E_{21})]_k\\
  =&  [\delta \lambda a_{11}E_{11}, \lambda^{-k}(E_{12} + E_{21})]_k
= \delta \lambda^{-k^2+1} a_{11}[E_{11}, E_{12} + E_{21}]_k.
\end{array}\eqno(2.5)$$
If $k$ is odd, one gets
$$\begin{array}{rl}
 &[E_{11}, E_{12} + E_{21}]_k\\
=& \sum_{i=0}^k(-1)^i{\rm C}_k^i(E_{12} + E_{21})^iE_{11}(E_{12} + E_{21})^{k-i}\\
=& \sum_{i=0}^{\frac{k-1}{2}}{\rm C}_k^{2i}E_{11}(E_{12} + E_{21}) - \sum_{i=0}^{\frac{k-1}{2}}{\rm C}_k^{2i+1} (E_{12} + E_{21})E_{11}\\
 =& \sum_{i=0}^{\frac{k-1}{2}}{\rm C}_k^{2i}E_{12} - \sum_{i=0}^{\frac{k-1}{2}}{\rm C}_k^{2i+1}  E_{21}\\
 =& 2^{k-1}(E_{12} - E_{21})\neq0;
\end{array}$$
If $k$ is even, one gets
$$\begin{array}{rl}
& [E_{11}, E_{12} + E_{21}]_k\\
=& \sum_{i=0}^k(-1)^i{\rm C}_k^i(E_{12} + E_{21})^iE_{11}(E_{12} + E_{21})^{k-i}\\
=&\sum_{i=0}^{\frac{k-1}{2}}{\rm C}_k^{2i} E_{11} - \sum_{i=0}^{\frac{k-1}{2}}{\rm C}_k^{2i+1} (E_{12} + E_{21})E_{11}(E_{12} + E_{21})\\
= &\sum_{i=0}^{\frac{k-1}{2}}{\rm C}_k^{2i} E_{11} - \sum_{i=0}^{\frac{k-1}{2}}{\rm C}_k^{2i+1}  E_{22}\\
 =& 2^{k-1}(E_{11} - E_{22})\neq0.
\end{array}$$
Hence, by Eq.(2.5), we see that $\delta = \lambda^{k^2-1}$ and then
$$\Phi(a_{11}E_{11}) =  \lambda^{k^2} a_{11}E_{11} + \mu_{a_{11}} I
$$ with $\mu_{a_{11}} = s_{22}$. Since
$\Phi(E_{11}) = \lambda E_{11} + \mu_1 I$ by Step 3, we get
$\lambda^{k^2-1} = 1$. Thus, we have
$$\Phi(a_{11}E_{11}) = \lambda a_{11}E_{11} + \mu_{a_{11}} I $$
for any $a_{11}\in{\mathbb F}$.

Next, we check that $\lambda^{k+1} = 1$.

Notice that
$$\begin{array}{rl}
&[E_{21}, E_{11} + E_{12}]_k \\ = &\sum_{i=0}^k(-1)^i{\rm C}_k^i(E_{11} + E_{12})^iE_{21}(E_{11} + E_{12})^{k-i} \\
= &\sum_{i=1}^{k-1}(-1)^i{\rm C}{_k^i}(E_{11} + E_{12})E_{21}(E_{11} + E_{12}) + E_{21} + E_{22} + (-1)^kE_{11}\\
= &\sum_{i=1}^{k-1}(-1)^i{\rm C}_k^i(E_{11} + E_{12}) + E_{21} + E_{22} + (-1)^kE_{11} \\
=&- E_{11} -(1+(-1)^k)E_{12}  + E_{21} + E_{22}.
\end{array}$$
By Step 4, one gets
$$\begin{array}{rl}
&[\Phi(E_{21}), \Phi(E_{11} + E_{12})]_k\\
= &[\lambda^{-k}E_{21}, \lambda E_{11} + \lambda^{-k}E_{12}]_k\\
=&\lambda^{-k}\sum_{i=0}^k(-1)^i{\rm C}_k^i(\lambda E_{11} + \lambda^{-k}E_{12})^iE_{21}(\lambda E_{11} + \lambda^{-k}E_{12})^{k-i}\\
=& \lambda^{-k}\sum_{i=1}^{k-1}(-1)^i{\rm C}_k^i(\lambda^i E_{11} +
\lambda^{i-1-k}E_{12})E_{21}(\lambda^{k-i} E_{11} +
\lambda^{-i-1}E_{12})\\
&+ \lambda^{-k}[E_{21}(\lambda^k E_{11} + \lambda^{-1}E_{12}) + (-1)^k(\lambda^k E_{11} + \lambda^{-1}E_{12})E_{21}]\\
=&\sum_{i=1}^{k-1}(-1)^i{\rm C}_k^i(\lambda^{-k-1} E_{11} +
\lambda^{-2-2k}E_{12}) + E_{21} + \lambda^{-k-1}E_{22} +
(-1)^k\lambda^{-k-1}E_{11}\\
=&-(1+(-1)^k)(\lambda^{-k-1} E_{11} + \lambda^{-2-2k}E_{12}) +
E_{21} + \lambda^{-k-1}E_{22} +
(-1)^k\lambda^{-k-1}E_{11}\\
=&- \lambda^{-k-1}E_{11} -(1+(-1)^k)\lambda^{-2k-2}E_{12} + E_{21} +
\lambda^{-k-1}E_{22}.
\end{array}$$
Since $[E_{21}, E_{11} + E_{12}]_k = [\Phi(E_{21}), \Phi(E_{11} +
E_{12})]_k$, comparing the above two equations gives
 $ \lambda^{k+1} = 1$.

\textbf{Step 6.} There exists a functional $h :{\mathcal
M}_2(\mathbb F)\to {\mathbb F}$ such that $\Phi(A) = \lambda A + h
(A)I$ holds for any $A\in{\mathcal M}_2(\mathbb F)$.

As $\lambda ^{-k}=\lambda$, by Step 4 and Step 5, we have that
$\Phi(a_{ij}E_{ij})= \lambda a_{ij}E_{ij} + \mu_{a_{ij}} I$ for any
$i, j\in\{1, 2\}$ and any $a_{ij}\in {\mathbb F}$. Then, for any
$A\in{\mathcal M}_2(\mathbb F)$,  writing $A = a_{11}E_{11} +
a_{12}E_{12} + a_{21}E_{21} + a_{22}E_{22}$ with $a_{ij}\in {\mathbb
F}$  and applying  Step 1,  there exists a scalar $c_A$ such that
$$\begin{array}{rl}
 \Phi(A) = & \Phi(a_{11}E_{11}) + \Phi(a_{12}E_{12}) + \Phi(a_{21}E_{21})+ \Phi(a_{22}E_{22}) + c_A I \\
=& \lambda a_{11}E_{11} + \mu_{a_{11}} I + \lambda a_{12}E_{12} +
\mu_{a_{12}} I + \lambda a_{21}E_{21} \\ &+ \mu_{a_{21}} I + \lambda
a_{22}E_{22} + \mu_{a_{22}} I + c_AI\\ = & \lambda A + (\mu_{a_{11}}
+ \mu_{a_{12}} + \mu_{a_{21}} + \mu_{a_{22}} + c_A) I = \lambda A +
\mu_A I.
\end{array}$$ For any $A \in{\mathcal M}_2(\mathbb F)$, let $h(A) = \mu_A$. Then $h
:{\mathcal M}_2(\mathbb F) \rightarrow {\mathbb F}$ is a functional
on ${\mathcal M}_2(\mathbb F)$ such that $\Phi(A) = \lambda A +
h(A)I$ holds for any $A\in{\mathcal
 M}_2(\mathbb F)$.

 The proof of Theorem 1 is completed. \hfill$\Box$


\end{document}